\documentclass{amsart} 
\usepackage{amssymb}
\usepackage{graphicx}
\usepackage{algpseudocode}
\usepackage[ruled]{algorithm}

\newtheorem{thm}{Theorem}[section] 
 \newtheorem{lem}[thm]{Lemma}
 \theoremstyle{remark}
\newtheorem{rem}[thm]{Remark} \newtheorem{ex}[thm]{Example}
 \theoremstyle{definition}
\newtheorem{defn}[thm]{Definition} \numberwithin{equation}{section}

\newcommand{\set}[1]{\left\{#1\right\}}

\renewcommand{\k}{\Bbbk}

\newcommand{\<}{\langle}
\renewcommand{\>}{\rangle}

\title[Algebraic Characterization of Uniquely Vertex Colorable Graphs]
 {Algebraic Characterization of \\ Uniquely Vertex Colorable Graphs}
 
 \author{Christopher J. Hillar}
\address{Department of Mathematics, Texas A\&M University, College Station, TX 77843}
\email{chillar@math.tamu.edu}

 \author{Troels Windfeldt}
\address{Department of Mathematical Sciences, University of Copenhagen, Denmark.}
\email{windfeldt@math.ku.dk}

 \thanks{The work of the first
  author was supported under a National Science Foundation Graduate
  Research Fellowship. The work of the second author was partially supported by Rejselegat for Matematikere. This work was conducted during the Special
  Semester on Gr\"obner Bases, February 1 -- July 31, 2006, organized by RICAM,
  Austrian Academy of Sciences, and RISC, Johannes Kepler University, Linz, Austria.}


\begin{document}

\begin{abstract}
The study of graph vertex colorability from an algebraic perspective has 
introduced novel techniques and algorithms into the field.  For instance, 
it is known that $k$-colorability of a graph $G$ is equivalent to the condition 
$1 \in I_{G,k}$ for a certain ideal $I_{G,k} \subseteq \k[x_1, \ldots, x_n]$.
In this paper, we extend this result 
by proving a general decomposition theorem for $I_{G,k}$.  
This theorem allows us to give an algebraic 
characterization of uniquely $k$-colorable graphs.
Our results also give algorithms for testing unique colorability.
As an application, we verify a counterexample to a conjecture
of Xu concerning uniquely $3$-colorable graphs without triangles.
\end{abstract}

\maketitle

\section{Introduction}

Let $G$ be a simple, undirected graph with vertices 
$V=\set{1,\ldots, n} $ and edges $ E $. The \emph{graph polynomial} of $G$
is given by
\[ f_G = \prod_{\substack{\set{i,j} \in E,\\i<j}} (x_i-x_j). \]
Fix a positive integer $k <n$, and let $C_k =\{c_1,\ldots,c_k\}$ be
a $k$-element set.  Each element of $C_k$ is called a \emph{color}. 
A (vertex) $k$-\textit{coloring} of $G$ is a map $\nu : V \to C_k$. We say that a
$k$-coloring $\nu$ is \emph{proper} if adjacent vertices receive
different colors; otherwise $\nu$ is called
\emph{improper}. The graph $G$ is said to be \emph{$k$-colorable} if
there exists a proper $k$-coloring of $G$.

Let $\k$ be an algebraically closed field of characteristic not
dividing $k$, so that it contains $k$ distinct $k$th roots of unity.  
Also, set $ R = \k[x_1, \ldots, x_n]$ 
to be the polynomial ring over $\k$ in indeterminates $x_1,\ldots,x_n$. 
Let $\mathcal{H} $ be the set of graphs with vertices $\{1,\ldots, n\}$ 
consisting of a clique of size $k+1$ and
isolated other vertices.
We will be interested in the following ideals of $R$:
\begin{equation*}\label{eqn_j}
\begin{split}
J_{n,k} = \ & \<f_H : H \in \mathcal{H}\>, \\
I_{n,k} = \ & \<x_i^k-1 : i \in V\>, \\
I_{G,k} = \ & I_{n,k}+\langle x_i^{k-1} + x_i^{k-2}x_j^{\phantom 2} + \cdots +  x_i^{\phantom 2}x_j^{k-2} + x_j^{k-1} :
\set{i,j} \in E \rangle. \\
\end{split}
\end{equation*}

One should think of (the zeroes of) $I_{n,k}$ and $I_{G,k}$ as representing 
$k$-colorings and proper $k$-colorings of the graph $G$, respectively
(see Section \ref{thmuc}).  The idea of using roots of unity and ideal theory
to study graph coloring problems seems to originate in Bayer's thesis \cite{bayer}, although
it has appeared in many other places, including the work of de Loera \cite{deloera}
and Lov\'asz \cite{lovasz}.  
These ideals are important because they allow for an algebraic formulation
of $k$-colorability.  The following theorem collects the results in the series of works
\cite{alon2,bayer,deloera, lovasz,mnuk}. 

\begin{thm}\label{thm_vc}
The following statements are equivalent:
\begin{enumerate}
\item The graph $ G $ is not $k$-colorable.
\item \text{\rm dim}$_{\k}R/I_{G,k} = 0$.
\item The constant polynomial $ 1 $ belongs to the ideal $ I_{G,k} $.
\item The graph polynomial $ f_G $ belongs to the ideal $ I_{n,k} $.
\item The graph polynomial $ f_G $ belongs to the ideal $ J_{n,k} $.
\end{enumerate}
\end{thm}

The equivalence between $(1)$ and $(3)$ is due to Bayer \cite[p. 109--112]{bayer}
(see also Chapter 2.7 of \cite{adams}).  Alon and Tarsi \cite{alon2} proved that
$(1)$ and $(4)$ are equivalent, but also de Loera \cite{deloera}
and Mnuk \cite{mnuk} have proved this using Gr\"obner basis methods.
The equivalence between $(1)$ and $(5)$ was proved by
Kleitman and Lov\'asz \cite{lovasz}.
We give a self-contained and simplified proof of Theorem \ref{thm_vc} in
Section \ref{algprelim}, in part to collect the many facts we need here.  

The next result says that the generators 
for the ideal $J_{n,k}$ in the above theorem are very
special.  A proof can be found in \cite{deloera}. 
(In Section \ref{algprelim}, we will review the relevant definitions 
regarding term orders and Gr\"obner bases). 

\begin{thm}[J. de Loera]\label{deloera}
The set of polynomials, $\{f_H : H \in \mathcal{H} \}$, is a universal Gr\"obner basis of $ J_{n,k} $.
\end{thm}

\begin{rem}
The set  $\mathcal{G} = \{x_1^k-1, \ldots, x_n^k-1\}$ is a universal Gr\"obner basis 
of $I_{n,k}$, but this follows easily since the leading terms of $\mathcal{G}$
are relatively prime, regardless of term order \cite[Theorem 1.7.4 and Lemma 3.3.1]{adams}.
\end{rem}

We say that a graph is \textit{uniquely $k$-colorable} if there is a unique proper $k$-coloring
up to permutation of the colors in $C_k$.  In this case,
partitions of the vertices into subsets having the same color
are the same for each of the $k!$ proper colorings of $G$.
A natural refinement of Theorem \ref{thm_vc}
would be an algebraic characterization of when a $k$-colorable
graph is uniquely $k$-colorable.  We provide such a 
characterization.  It will be a corollary to our main theorem (Theorem \ref{bigmainthm})
that decomposes the ideal $I_{G,k}$ into an intersection of simpler ``coloring ideals". 
To state the theorem, however, we need to introduce some notation.

Let $\nu$ be a proper $k$-coloring of a graph $G$.
Also, let $l \leq k$ be the number of distinct colors in $\nu(V)$.  
The \textit{color class} $cl(i)$ of a vertex $i \in V$ is the set of vertices with the 
same color as $i$, and the \textit{maximum} of a color class is the largest
vertex contained in it.  We set $m_1 < m_2 < \cdots < m_l = n$ to be the maximums
of the $l$ color classes.

For a subset $U \subseteq V$ of the vertices, let $h_U^d$ be the sum of all monomials 
of degree $d$ in the indeterminates $\{x_i : i\in U \}$.  We also set $h_U^0 = 1$.

\begin{defn}[$\nu$-bases]\label{nubasis}
Let $\nu$ be a proper $k$-coloring of a graph $G$.
For each vertex $i \in V$, define a polynomial $g_i$ as follows:
\begin{equation}\label{gidefn}
g_i  = \begin{cases} 
\ x_i^k - 1   &  \ \ \text{if $i = m_l$}, \\
\  h_{\{m_j,\ldots,m_l\}}^{k-l+j}   & \ \ \text{if $i = m_j$ for some $j \neq l$},  \\
\ x_i - x_{\max{cl(i)}} &  \ \ \text{otherwise}.
 \end{cases}
\end{equation}
The collection $\{g_1,\ldots,g_n\}$ is called a 
$\nu$-\textit{basis} for the graph $G$ with respect to the proper coloring $\nu$.
\end{defn}

As we shall soon see, this set is a (minimal) Gr\"obner basis; its initial ideal
is generated by the relatively prime monomials 
\[\{x_{m_1}^{k-l+1},x_{m_2}^{k-l+2},\ldots,x_{m_l}^k\} \ \text{ and } \ \{x_i : i \neq m_j \text{ for any } j\}.\]
A concrete instance of this construction may be found in Example \ref{decompex} below.

\begin{rem}It is easy to see that the map $\nu \mapsto \{g_1,\ldots,g_n\}$  
depends only on  how $\nu$ partitions $V$ into color classes $cl(i)$.  
In particular, if $G$ is uniquely $k$-colorable, then there is a unique
such set of polynomials $\{g_1,\ldots,g_n\}$ that corresponds to $G$.
\end{rem}

This discussion prepares us to make the following definition.

\begin{defn}[Coloring Ideals]\label{colorideal}
Let $\nu$ be a proper $k$-coloring of a graph $G$.
The \textit{k-coloring ideal} (or simply \textit{coloring ideal}
if $k$ is clear from the context) associated to $\nu$ is the ideal
\[A_{\nu} = \<g_1,\ldots,g_n\>,\] where the $g_i$ are given by 
(\ref{gidefn}).
\end{defn}


In a precise way to be made clear later (see Lemma \ref{giAlem}), 
the coloring ideal associated to $\nu$ algebraically encodes the proper $k$-coloring of $G$ by $\nu$
(up to relabeling of the colors).
We may now state our main theorem.  

\begin{thm}\label{bigmainthm}
Let $G$ be a simple graph with $n$ vertices.  Then
\[ I_{G,k} = \bigcap_{\nu} A_{\nu},\]
where $\nu$ runs over all proper $k$-colorings of $G$.
\end{thm}

\begin{ex}\label{decompex}
Let $G = (\{1, 2, 3\}, \{\{1,2\}, \{2,3\}\}) $ be the path graph on three vertices, and let $k=3$. There are essentially two proper $3$-colorings of $G$: the one where vertices 1 and 3 receive the same color, and the one where all the vertices receive different colors. If we denote by $\nu_1$ the former, and by $\nu_2$ the latter, then according to Definition \ref{colorideal}, we have:
\begin{equation*}
\begin{split}
A_{\nu_1} = \ & \<x_3^3-1, x_2^2+x_2 x_3+x_3^2, x_1-x_3\>, \\
A_{\nu_2} = \ & \<x_3^3-1, x_2^2+x_2 x_3+x_3^2, x_1+x_2+x_3\>. \\
\end{split}
\end{equation*}
The intersection $A_{\nu_1} \cap A_{\nu_2}$ is equal to the graph ideal,
\[
I_{G,3} = \< x_1^3-1, x_2^3-1, x_3^3-1, x_1^2+x_1 x_2+x_2^2, x_2^2+x_2 x_3+x_3^2 \>,
\]
as predicted by Theorem \ref{bigmainthm}.
\end{ex}

Two interesting special cases of this theorem are the following.
When $G$ has no proper $k$-colorings, Theorem \ref{bigmainthm}
says that $I_{G,k} = \<1\>$  in accordance with Theorem \ref{thm_vc}.
And for a graph that is uniquely $k$-colorable, all of the ideals $A_{\nu}$ are the
same.  This observation allows us to use Theorem \ref{bigmainthm} 
to give the following algebraic characterization of uniquely colorable graphs.

\begin{thm}\label{thm_uvc}
Suppose $\nu$ is a $k$-coloring of $G$ that uses all $k$ colors,
and let $g_1, \ldots, g_n$ be given by (\ref{gidefn}). Then the following statements are equivalent:
\begin{enumerate}
\item The graph $G$ is uniquely $k$-colorable.
\item The polynomials $g_1, \ldots, g_n$ generate the ideal $ I_{G,k} $.
\item The polynomials $g_1, \ldots, g_n$ belong to the ideal $ I_{G,k} $.
\item The graph polynomial $f_G$ belongs to the ideal $I_{n,k}:\<g_1, \ldots, g_n\>$.
\item \text{\rm dim}$_{\k}R/I_{G,k} = k!$.
\end{enumerate}
\end{thm}

There is also a partial analogue to Theorem \ref{deloera} that refines Theorem \ref{thm_uvc}.
This result gives us an algorithm for determining unique $k$-colorability
that is independent of the knowledge of a proper coloring.  To state it, we need only make a slight 
modification of the polynomials in (\ref{gidefn}).  Suppose that $\nu$ is a proper
coloring with $l = k$ (for instance, this holds when $G$ is uniquely $k$-colorable).  
Then, for $i \in V$ we define:
\begin{equation}\label{gidefnnunique}
\tilde{g}_i  = \begin{cases} 
\ x_i^k - 1   &  \ \ \text{if $i = m_l$}, \\
\  h_{\{m_j,\ldots,m_l\}}^{j}   & \ \ \text{if $i = m_j$ for some $j \neq l$},  \\
 \  h_{\{i,m_2,\ldots,m_l\}}^{1}   & \ \ \text{if $i \in cl(m_1)$},  \\
\ x_i - x_{\max{cl(i)}} &  \ \ \text{otherwise}.
 \end{cases}
\end{equation}
We call the set $\{\tilde{g}_1,\ldots,\tilde{g}_n\}$ a \textit{reduced} $\nu$-\textit{basis}.

\begin{rem}\label{nubasisrem}
When $l = k$, the ideals generated by the polynomials in (\ref{gidefn}) 
and in  (\ref{gidefnnunique}) are the same.  This follows because for $i \in cl(m_1) \setminus \{m_1\}$,
we have $\tilde{g}_i - \tilde{g}_{m_1} = x_i - x_{m_1} = g_i$.
\end{rem} 

\begin{thm}\label{thmGB}
A graph $G$ with $n$ vertices is uniquely $k$-colorable if and only if the reduced Gr\"obner basis 
for $ I_{G,k} $ with respect to any term order with $x_n \prec \cdots \prec x_1$
has the form $\{\tilde{g}_1,\ldots, \tilde{g}_n\}$ for polynomials as in $(\ref{gidefnnunique})$.
\end{thm}

\begin{rem}
It is not difficult to test whether a Gr\"obner basis is of the form given
by $(\ref{gidefnnunique})$.  Moreover, the unique coloring can be easily recovered
from the reduced Gr\"obner basis.
\end{rem}

In Section \ref{sec_alg}, we shall discuss the tractability of our algorithms.
We hope that they might be used to perform experiments for raising and settling
problems in the theory of (unique) colorability.

\begin{ex}\label{firstex}
We present an example of a uniquely $3$-colorable graph on
$n = 12$ vertices and give the polynomials $\tilde{g}_1, \ldots, \tilde{g}_n$ from 
Theorem \ref{thmGB}.

\begin{figure}[!htbp]
\begin{center}
\includegraphics[scale=1.0]{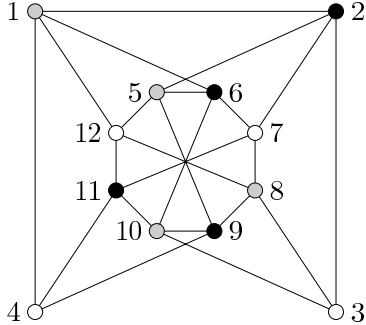}
\end{center}
\caption{A uniquely $3$-colorable graph \cite{chao}.}  
\label{fig.unique3colorex}
\end{figure}

Let $G$ be the graph given in Figure \ref{fig.unique3colorex}.  The 
indicated $3$-coloring partitions $V$ into $k = l = 3$ color classes with 
$(m_1,m_2,m_3) = (10,11,12)$.
The following set of $12$ polynomials is the
reduced Gr\"obner basis for the ideal $I_{G,k}$ with respect to any
term ordering with $x_{12} \prec \cdots \prec x_1$.  The leading
terms of each $\tilde{g}_i$ are underlined.
\[
\begin{split}
\{
& \underline{x_{12}^3}-1\text{, }
\underline{x_{7}}-x_{12}\text{, }
\underline{x_{4}}-x_{12}\text{, }
\underline{x_{3}}-x_{12}\text{, } \\
& \underline{x_{11}^2}+x_{11}x_{12}+x_{12}^2\text{, }
\underline{x_{9}}-x_{11}\text{, } 
\underline{ x_{6}}-x_{11}\text{, }
 \underline{x_{2}}-x_{11}\text{, }\\
& \underline{x_{10}}+x_{11}+x_{12}\text{, }
\underline{x_{8}}+x_{11}+x_{12}\text{, }
\underline{x_{5}}+x_{11}+x_{12}\text{, }
\underline{x_1}+x_{11}+x_{12}
\}.
\end{split}
\]
Notice that the leading terms of the polynomials in each line above 
correspond to the different color classes of this coloring of $G$.  \qed
\end{ex}

The organization of this paper is as follows.
In Section \ref{algprelim}, we discuss some of the algebraic tools 
that will go into the proofs of our main results.  Section \ref{thmuc} is devoted
to a proof of Theorem \ref{thm_vc}, and in Sections \ref{sec_decomp} and \ref{sec_uvc}, 
we present proofs for Theorems \ref{bigmainthm}, \ref{thm_uvc}, and \ref{thmGB}.  
Theorems \ref{thm_vc} and \ref{thm_uvc} give algorithms for testing 
$k$-colorability and unique $k$-colorability of graphs, and we
discuss the implementation of them in Section \ref{sec_alg}, along
with a verification of a counterexample \cite{akbari1} 
to a conjecture \cite{chao,DA,xu} by Xu concerning uniquely 3-colorable graphs without triangles.

\section{Algebraic Preliminaries}\label{algprelim}

We briefly review the basic concepts of commutative algebra
that will be useful for us here.  We refer to \cite{Cox2} or \cite{Cox1} for more details.
Let $I$ be an ideal of $R = \k[x_1, \ldots, x_n]$.  The \textit{variety} $V(I)$ of $I$ is the set of 
points in $\k^n$ that are zeroes of all the polynomials in $I$.  Conversely, the \textit{vanishing
ideal} $I(V)$ of a set $V \subseteq \k^n$ is the ideal of those polynomials vanishing on all 
of $V$.  These two definitions are related by way of $V(I(V)) = V$ and $I(V(I)) = \sqrt I$, in which 
\[\sqrt I = \{ f : f^n \in I \text{ for some } n\}\] is the \textit{radical} of $I$.  
The ideal $I$ is said to be of \textit{Krull dimension zero} (or simply
\textit{zero-dimensional}) if $V(I)$ is finite.
A \textit{term order} $\prec$ for the monomials of $R$ is a well-ordering which is multiplicative
($u \prec v \Rightarrow wu \prec wv$ for monomials $u,v,w$) 
and for which the constant monomial $1$ is smallest.  The \textit{initial term} (or \textit{leading monomial}) 
$in_{\prec}(f)$ of a polynomial $f \in R$ is the largest monomial in $f$ with
respect to $\prec$.  The \textit{standard monomials} $\mathcal{B}_{\prec}(I)$
of $I$ are those monomials which are not the leading monomials of any polynomial in $I$.

Many arguments in commutative algebra and algebraic geometry are simplified when 
restricted to radical, zero-dimensional ideals (resp. multiplicity-free,
finite varieties), and those found in this paper are not exceptions.  The 
following fact is useful in this regard.

\begin{lem}\label{radlem}
Let $I$ be a zero-dimensional ideal and fix a term order $\prec$.  
Then $\dim_{\k} R/I = |\mathcal{B}_{\prec}(I)| \geq |V(I)|$.  Furthermore, the following are equivalent:
\begin{enumerate}
\item  $I$ is a radical ideal (i.e., $I = \sqrt{I}$).
\item  $I$ contains a univariate square-free polynomial in each indeterminate.
\item  $|\mathcal{B}_{\prec}(I)| = |V(I)|$.
\end{enumerate}
\end{lem}
\begin{proof}
See \cite[p. 229, Proposition 4]{Cox1} and \cite[pp. 39--41, Proposition 2.7 and Theorem 2.10]{Cox2}.
\end{proof}
A finite subset $\mathcal{G}$ of an ideal $I$ is a \textit{Gr\"obner basis} 
(with respect to $\prec$) if the \textit{initial ideal}, 
\[ in_{\prec}(I) = \< in_{\prec}(f) : f \in I  \>, \] is
generated by the initial terms of elements of $\mathcal{G}$.
It is called \textit{minimal} if no leading term of $f \in G$
divides any other leading term of polynomials in $G$.
Furthermore, a \textit{universal Gr\"obner basis} is a set of polynomials 
which is a Gr\"obner basis with respect to all term orders.
Many of the properties of $I$ and $V(I)$ can be calculated by finding 
a Gr\"obner basis for $I$, and such generating sets are fundamental
for computation (including the algorithms presented in the last section).

Finally, a useful operation on two ideals $I$ and $J$ is the construction 
of the \textit{colon ideal} $I:J = \{h \in R : hJ \subseteq I \}$.  If $V$ and
$W$ are two varieties, then the colon ideal 
\begin{equation}\label{colonidl}
I(V) : I(W) = I(V \backslash W)
 \end{equation}
corresponds to a set difference \cite[p. 193, Corollary 8]{Cox1}.

\section{Vertex Colorability}\label{thmuc}

In what follows, the set of colors $C_k$ will be the set of $k$th roots of unity,
and we shall freely speak of points in $\k^n$ with all coordinates in $C_k$
as colorings of $G$.  In this case, a point $(v_1,\ldots,v_n) \in \k^n$ corresponds to
a coloring of vertex $i$ with color $v_i$ for $i = 1,\ldots,n$.
The varieties corresponding to the ideals $I_{n,k}$, $I_{G,k}$, and $I_{n,k}+\<f_G\>$
partition the $k$-colorings of $G$ as follows.

\begin{lem}\label{varietycolcorres}
The varieties $V(I_{n,k})$, $V(I_{G,k})$, and $V(I_{n,k}+\<f_G\>)$
are in bijection with all, the proper, and the improper $k$-colorings of $G$, respectively.
\end{lem}
\begin{proof}
The points in $V(I_{n,k})$ are all $n$-tuples of $k$th roots of unity and therefore 
naturally correspond to all $k$-colorings of $G$.  Let $\mathbf{v} = (v_1,\ldots,v_n) \in V(I_{G,k})$; we must show that 
it corresponds to a proper coloring of $G$.  Let $\{i,j\} \in E$ and set 
\[q_{ij} = \frac{x_i^k-x_j^k}{x_i-x_j} \in I_{G,k}.\]
If $v_i = v_j$, then $q_{ij}(\mathbf{v}) = k v_i^{k-1} \neq 0$.  Thus, the coloring
$\mathbf{v}$ is proper.
Conversely, suppose that $\mathbf{v} = (v_1,\ldots,v_n)$ is a proper coloring of $G$.  Then, since
\[q_{ij}(\mathbf{v})(v_i-v_j) =  (v_i^k-v_j^k) =  1-1 = 0,\]
it follows that for $\{i,j\} \in E$, we have $q_{ij}(\mathbf{v}) = 0$.
This shows that $\mathbf{v} \in V(I_{G,k})$.
If $\mathbf{v}$ is an improper coloring, then it is easy to see that
$f_G(\mathbf{v}) = 0$.   Moreover, any $\mathbf{v} \in V(I_{n,k})$ for which
$f_G(\mathbf{v}) = 0$ has two coordinates, corresponding to an edge in $G$, that are equal.
\end{proof}

The next result follows directly from Lemma \ref{radlem}.  It will prove useful
in simplifying many of the proofs in this paper.

\begin{lem}\label{radIcor}
The ideals $I_{n,k}$, $I_{G,k}$, and $I_{n,k}+\<f_G\>$ are radical.
\end{lem}

We next describe a relationship
between $I_{n,k}$, $I_{G,k}$, and  $I_{n,k} + \<f_G\>$.  

\begin{lem}\label{lem_duality}
$I_{n,k}:I_{G,k}=I_{n,k}+\<f_G\>$.
\end{lem}
\begin{proof}
Let $V$ and $W$ be the set of all colorings and proper colorings, respectively,
of the graph $G$.  Now apply Lemma \ref{varietycolcorres} and Lemma \ref{radIcor} to
equation (\ref{colonidl}).
\end{proof}

The vector space dimensions of the residue rings corresponding to these ideals 
are readily computed from the above discussion.  Recall that
the \textit{chromatic polynomial} $\chi_G$ is the univariate polynomial
for which $\chi_G(k)$ is the number of proper $k$-colorings of $G$. 

\begin{lem}\label{lem_chi}
Let $ \chi_G $ be the chromatic polynomial of $ G $. Then
\begin{equation*}\label{eqn_good}
\hspace*{-0.58cm}\chi_G(k) = \dim_{\k}R / I_{G,k},
\end{equation*}
\begin{equation*}\label{eqn_bad}
k^n-\chi_G(k) = \dim_{\k}R /(I_{n,k}+\<f_G\>).
\end{equation*}
\end{lem}
\begin{proof}
Both equalities follow from Lemmas \ref{radlem} and \ref{varietycolcorres}.
\end{proof}

Let $K_{n,k}$ be the ideal of all polynomials $f \in R$ such that
$f(v_1,\ldots,v_n) = 0$ for any $(v_1,\ldots,v_n) \in \k^n$ with
at most $k$ of the $v_i$ distinct.  Clearly, $J_{n,k} \subseteq K_{n,k}$.
We will need the following result of Kleitman  and Lov\'asz \cite{lovasz}.

\begin{thm}[Kleitman-Lov\'asz]\label{lovkleit}
The ideals $K_{n,k}$ and $J_{n,k}$ are the same.
\end{thm}

We now prove Theorem \ref{thm_vc}.  We feel that it
is the most efficient proof of this result.

\begin{proof}[Proof of Theorem \ref{thm_vc}]\mbox{}
$(1) \Rightarrow (2) \Rightarrow (3)$: Suppose that $G$ is not $k$-colorable. 
Then it follows from Lemma \ref{lem_chi} that 
$\dim_{\k}R / I_{G,k} = 0$ and so $1 \in I_{G,k}$.  

$(3) \Rightarrow (4)$:
Suppose that $I_{G,k}  = \<1\>$ so that $I_{n,k}:I_{G,k} = I_{n,k}$. Then
Lemma \ref{lem_duality} implies that \mbox{$I_{n,k} + \<f_G\>=I_{n,k}$} and hence $f_G\in I_{n,k}$.
 
 $(4) \Rightarrow (1)$:  Assume that $f_G $ 
belongs to the ideal $I_{n,k}$. Then $I_{n,k} + \<f_G\> = I_{n,k}$, and it follows from 
Lemma \ref{lem_chi} that $k^n-\chi_G(k)=k^n$.  Therefore, $\chi_G(k)=0$ as desired.

$(5) \Rightarrow (1)$:  Suppose that $f_G \in J_{n,k}$.  Then from Theorem \ref{lovkleit},
there can be no proper coloring $\mathbf{v}$ (there are at most $k$ distinct coordinates).

$(1) \Rightarrow (5)$:  If $G$ is not $k$-colorable, then for every substitution
$\mathbf{v} \in \k^n$ with at most $k$ distinct coordinates, we must
have $f_G(\mathbf{v}) = 0$.  It follows that $f_G \in J_{n,k}$ from Theorem \ref{lovkleit}.
\end{proof}

\section{Coloring Ideals}
\label{sec_decomp}

In this section, we study the $k$-coloring ideals $A_{\nu}$ mentioned in the introduction
and prove Theorem \ref{bigmainthm}.
Let $G$ be a graph with proper coloring $\nu$, and let $l \leq k$ be the 
number of distinct colors in $\nu(V)$.  For each vertex $i \in V$, we assign
polynomials $g_i$ and $\tilde{g}_i$ as in equations (\ref{gidefn}) and (\ref{gidefnnunique}).  
One should think (loosely) of the first case of (\ref{gidefn}) as corresponding 
to a choice of a color for the last vertex;  the second, to subsets of vertices in different 
color classes; and the third, to the fact that elements in the same color class should have the same color. These polynomials encode the coloring $\nu$ algebraically in a
computationally useful way (see Lemmas \ref{GBgilem} and \ref{giAlem} below).
We begin by showing that the 
polynomials $g_i$ are a special generating set for the coloring ideal $A_{\nu}$.

Recall that a \textit{reduced Gr\"obner basis} $\mathcal{G}$ is
a Gr\"obner basis such that $(1)$ the coefficient of $in_{\prec}(g)$ for each $g\in \mathcal{G}$ 
is $1$ and $(2)$ the leading monomial of any $g \in \mathcal{G}$ does not divide any 
monomial occurring in another polynomial in $\mathcal{G}$.  Given a term order,
reduced Gr\"obner bases exist and are unique.  

\begin{lem}\label{GBgilem}
Let $\prec$ be any term order with $x_n \prec  \cdots \prec x_1$. Then
the set of polynomials $\{g_1,\ldots,g_n\}$ is a minimal Gr\"obner basis 
with respect to $\prec$ for the ideal $A_{\nu} = \<g_1,\ldots,g_n\>$  it generates.
Moreover, for this ordering, the set $\{\tilde{g}_1,\ldots,\tilde{g}_n\}$ is a reduced 
Gr\"obner basis for $\<\tilde{g}_1,\ldots,\tilde{g}_n\>$.
\end{lem}
\begin{proof}
Since the initial term of each $g_i$ (resp. $\tilde{g}_i$) is a power of $x_i$, each pair of leading
terms is relatively prime.  It follows that these polynomials form
a Gr\"obner basis for the ideal they generate. 
By inspection, it is easy to see that the set of polynomials given by (\ref{gidefn})
(resp. (\ref{gidefnnunique})) is minimal (resp. reduced).
\end{proof}

The following innocuous-looking fact is a very important ingredient in the 
proof of Lemma  \ref{giAlem}.

\begin{lem}\label{hlemma}
Let $U$ be a subset of $\{1,\ldots,n\}$, and suppose that $\{i,j\}\subseteq U$. Then
\begin{equation}
(x_i-x_j)h_U^d = h_{U \backslash \{j\}}^{d+1}-h_{U\backslash \{i\}}^{d+1},
\end{equation}
for all nonnegative integers $d$.
\end{lem}
\begin{proof}
The first step is to note that the polynomial
\[
x_i h_U^d + h_{U\backslash\set{i}}^{d+1}
\]
is symmetric in the indeterminants $\{x_\ell : \ell\in U \}$. This 
follows from the polynomial identity
\[
h_U^{d+1} - h_{U\backslash\set{i}}^{d+1} = x_i h_U^d,
\]
and the fact that $h_U^{d+1}$ is symmetric in the indeterminants $\{x_\ell :
\ell\in U \}$. Let $\sigma$ be the permutation $(i \ j)$, and
notice that
\[
x_i h_U^d + h_{U\backslash\set{i}}^{d+1} =
\sigma\left(x_i h_U^d + h_{U\backslash\set{i}}^{d+1}\right) =
x_j h_U^d + h_{U\backslash\set{j}}^{d+1}.
\]
This completes the proof.
\end{proof}

We shall also need the following fact that gives explicit 
representations of some of the generators of $I_{n,k}$
in terms of those of $A_{\nu}$.

\begin{lem}\label{Inklem}
For each $i = 1,\ldots,l$, we have
\begin{equation}\label{telescopesum}
x_{m_i}^k - 1 = x_n^k-1 + \sum_{t=i}^{l-1} \left[ \prod_{j = t+1}^l 
\left(x_{m_i}- x_{m_j} \right) \right] h_{\{m_{t},\ldots,m_l\}}^{k-l+t}.
\end{equation}
\end{lem}
\begin{proof}
To verify (\ref{telescopesum}) for a fixed $i$, we will use Lemma \ref{hlemma} and induction to 
prove that for each positive integer $s \leq l-i$, the sum on the  right hand-side above is equal to
\begin{equation}\label{telescopeinduction}
\prod_{j = s+i}^l \left(x_{m_i}- x_{m_j} \right) h_{\{m_i,m_{s+i},\ldots,m_l\}}^{k-l+s+i-1} 
+ \sum_{t=s+i}^{l-1} \left[ \prod_{j = t+1}^l
\left(x_{m_i}- x_{m_j} \right) \right] h_{\{m_{t},\ldots,m_l\}}^{k-l+t}.
\end{equation}
For $s=1$, this is clear as (\ref{telescopeinduction}) is exactly the
sum on the right-hand side of (\ref{telescopesum}).  In general, using Lemma \ref{hlemma}, the first
term on the left hand side of (\ref{telescopeinduction}) is
\[   \prod_{j = s+1+i}^l \left(x_{m_i}- x_{m_j} \right) \left(h_{\{m_i,m_{s+1+i},\ldots,m_l\}}^{k-l+s+i} -
h_{\{m_{s+i},\ldots,m_l\}}^{k-l+s+i} \right),\]
which is easily seen to cancel the first summand in the 
sum found in (\ref{telescopeinduction}).  

Now, equation
(\ref{telescopeinduction}) with $s = l-i$ gives us that the right hand side
of (\ref{telescopesum}) is 
\[ x_n^k-1 + (x_{m_i}- x_{m_l} ) h_{\{m_i,m_l\}}^{k-1} = x_n^k-1 + x_{m_i}^k - x_n^k = x_{m_i}^k-1,\]
proving the claim (recall that $m_l = n$).
\end{proof}

That the polynomials $g_1,\ldots,g_n$ represent an algebraic encoding of the coloring $\nu$
is explained by the following technical lemma.

\begin{lem}\label{giAlem}
Let $g_1,\ldots,g_n$ be given as in (\ref{gidefn}).  Then the following 
three properties hold for the ideal $A_{\nu} = \<g_1,\ldots,g_n\>$:
\begin{enumerate}
\item $I_{G,k} \subseteq  A_{\nu}$,
\item $A_{\nu}$ is radical, 
\item $|V(A_{\nu})| = \prod_{j=1}^{l} (k-l+j)$.
\end{enumerate}
\end{lem}
\begin{proof}
First assume that $I_{G,k} \subseteq A_{\nu}$.  Then $A_{\nu}$ is radical from
Lemma \ref{radlem}, and the number of standard monomials of $A_{\nu}$ 
(with respect to any ordering $\prec$ as in Lemma \ref{GBgilem}) is equal to $|V(A_{\nu})|$.
Since $\{g_1,\ldots,g_n\}$ is a Gr\"obner basis for $A_{\nu}$ and the initial ideal is
generated by the monomials \[\{x_{m_1}^{k-l+1},x_{m_2}^{k-l+2},\ldots,x_{m_l}^k\} \ \text{ and } \ \{x_i : i \neq m_j \text{ for any } j\},\]
it follows that $|\mathcal{B}_{\prec}(A_{\nu})| = \prod_{j=1}^{l} (k-l+j)$. This proves $(3)$.  

We now prove statement $(1)$.   From Lemma \ref{Inklem}, it follows that 
$x_i^k-1 \in A$ when $i \in \{m_1,\ldots,m_l\}$.
It remains to show that $x_i^k-1 \in A_{\nu}$ for all vertices not in $\{m_1,\ldots,m_l\}$.
Let $f_i = x_i - x_{\max cl(i)}$ and
notice that \[ x_{\max cl(i)}^k - 1 = (x_i-f_i)^k - 1 = x_i^k - 1 + f_i h \in A_{\nu}\]
for some polynomial $h$.  It follows that $x_i^k - 1 \in A_{\nu}$.

Finally, we must verify that the other generators of $I_{G,k}$ are in $A_{\nu}$.
To accomplish this, we will prove the following stronger statement:
\begin{equation}\label{subsetinA}
U \subseteq \{m_1,\ldots,m_l\} \text{ with } |U| \geq 2 \  \Longrightarrow \  h_{U}^{k+1-|U|} \in A_{\nu}.
\end{equation}
We downward induct on $s = |U|$.
In the case $s = l$, we have $U =\{m_1,\ldots,m_l\}$.
But then as is easily checked $g_{m_1} = h_{U}^{k+1-|U|} \in A_{\nu}$.
For the general case, we will show that if one polynomial $h_{U}^{k+1-|U|}$ is in $A_{\nu}$,
with $|U| = s < l$, then $h_{U}^{k+1-|U|} \in A_{\nu}$ for any subset $U \subseteq \{m_1,\ldots,m_l\}$ of
cardinality $s$.  In this regard, suppose that $h_{U}^{k+1-|U|} \in A_{\nu}$ for
a subset $U$ with $|U| = s < l$.  Let $u \in U$ and 
$v \in  \{m_1,\ldots,m_l\} \backslash U$, and examine the following equality
(using Lemma \ref{hlemma}):
\[ (x_u - x_v) h_{\{v\} \cup U}^{k-s} = h_{U}^{k-s+1} - h_{\{v\} \cup U \backslash \{u\} }^{k-s+1}.\]
By induction, the left hand side of this equation is in $A_{\nu}$ and therefore
the assumption on $U$ implies that
\[ h_{\{v\} \cup U \backslash \{u\}}^{k-s+1}  \in A_{\nu}.\]
This shows that we may replace any element of $U$ with
any element of $\{m_1,\ldots,m_l\}$.  Since there is a subset 
$U$ of size $s$ with  $h_{U}^{k+1-|U|} \in A_{\nu}$ (see (\ref{gidefn})), 
it follows from this that we have $h_U^{k+1-|U|} \in A_{\nu}$ for any subset $U$ of size $s$.
This completes the induction.

A similar trick as before using polynomials $x_i - x_{\max cl(i)} \in A_{\nu}$
proves that we may replace in (\ref{subsetinA})
the requirement that $U \subseteq \{m_1,\ldots,m_l\}$ with one
that says that $U$ consists of vertices in different color classes.
If $\{i,j\} \in E$, then $i$ and $j$ are in different color classes,
and therefore the generator $h_{\{i,j\}}^{k-1} \in I_{G,k}$ is in $A_{\nu}$.
This finishes the proof of the lemma.
\end{proof}

\begin{rem}\label{Alemremark}
Property $(1)$ in the lemma says that $V(A_{\nu})$ contains only proper colorings 
of $G$ while properties $(2)$ and $(3)$ say that, up to relabeling the
colors, the zeroes of the polynomials $g_1,\ldots,g_n$ correspond
to the single proper coloring given by $\nu$.  The lemma also 
implies that the polynomials $\{g_1,\ldots,g_n\}$ form a complete intersection.
\end{rem}

The decomposition theorem for $I_{G,k}$ mentioned in the introduction
now follows easily from the results of this section.

\begin{proof}[Proof of Theorem \ref{bigmainthm}]
By Lemmas \ref{varietycolcorres} and \ref{giAlem}, we have
\begin{equation*}
\begin{split}
V(I_{G,k}) = \ & \bigcup_{\nu} V(A_{\nu}),  \\
\end{split}
\end{equation*}
where $\nu$ runs over all proper $k$-colorings of $G$.  
Since the ideals $I_{G,k}$ and $A_{\nu}$ are radical by 
Lemmas \ref{radIcor} and \ref{giAlem}, it follows that:
\begin{equation*}
\begin{split}
I_{G,k} = \ & I(V(I_{G,k})) \\
= \ & I\bigcup_{\nu} V(A_{\nu}) \\
= \ & \bigcap_{\nu} I(V(A_{\nu})) \\
= \ & \bigcap_{\nu} A_{\nu}.  \\
\end{split}
\end{equation*}
This completes the proof.
\end{proof}

\section{Unique Vertex Colorability}
\label{sec_uvc}

We are now in a position to prove our characterizations of uniquely
$k$-colorable graphs.

\begin{proof}[Proof of Theorem \ref{thm_uvc}]\mbox{}
$(1) \Rightarrow (2) \Rightarrow (3)$: Suppose the graph $G$ is uniquely $k$-colorable and 
construct the set of $g_i$ from (\ref{gidefn}) using the proper $k$-coloring $\nu$.  
By Theorem \ref{bigmainthm}, it follows that $I_{G,k} = A_{\nu}$, and thus
the $g_i$ generate $I_{G,k}$.



$(3) \Rightarrow (4)$: Suppose that $A_{\nu} = \<g_1,\ldots,g_n\> \subseteq I_{G,k}$.  From
Lemma \ref{lem_duality}, we have 
\begin{equation*}
\begin{split}
I_{n,k} + \<f_G\> =   I_{n,k} : I_{G,k}  \subseteq   I_{n,k} : A_{\nu}.\\
\end{split}
\end{equation*} 
This proves that $f_G \in I_{n,k}:\<g_1, \ldots, g_n\>$.

$(4) \Rightarrow (5) \Rightarrow (1)$:  Assume that $f_G \in I_{n,k}:\<g_1, \ldots, g_n\>$.  Then, 
\[I_{n,k} : I_{G,k} =  I_{n,k} + \<f_G\> \subseteq  I_{n,k}:\<g_1, \ldots, g_n\>.\]
Applying Lemmas \ref{colonidl} and \ref{giAlem}, we have
\begin{equation}\label{inequvc}
 k^n - k! = |V(I_{n,k}) \backslash V(A_{\nu})|  = |V(I_{n,k} : A_{\nu})| \leq |V(I_{n,k} : I_{G,k})| \leq k^n-k!,
\end{equation}
since the number of improper colorings is at most $k^n-k!$.
It follows that equality holds throughout (\ref{inequvc})
so that the number of proper colorings is $k!$. Therefore,
we have \text{\rm dim}$_{\k}R/I_{G,k} = k!$ from Lemma \ref{lem_chi} 
and $G$ is uniquely $k$-colorable.
\end{proof}

\begin{proof}[Proof of Theorem \ref{thmGB}]
Suppose that the reduced Gr\"obner basis of $I_{G,k}$ with respect to a
term order with $x_n \prec \cdots \prec x_1$ has the form 
$\{\tilde{g}_1,\ldots,\tilde{g}_n\}$ as in $(\ref{gidefnnunique})$.  Also, let 
$\{g_1,\ldots,g_n\}$ be the $\nu$-basis (\ref{gidefn}) corresponding to the $k$-coloring
$\nu$ read off from $\{\tilde{g}_1,\ldots,\tilde{g}_n\}$.
By Remark \ref{nubasisrem}, we have  $\<g_1,\ldots,g_n\> = \<\tilde{g}_1,\ldots,\tilde{g}_n\>$.
It follows that $G$ is uniquely $k$-colorable from $(2) \Rightarrow (1)$ of Theorem \ref{thm_uvc}.   
For the other implication, by Lemma \ref{GBgilem},
it is enough to show that  $A_{\nu} = \<g_1,\ldots,g_n\> = I_{G,k}$, which is  $(1) \Rightarrow (2)$ 
in Theorem \ref{thm_uvc}.
\end{proof}

\section{Algorithms and Xu's Conjecture}
\label{sec_alg}

In this section we describe the algorithms implied by Theorems \ref{thm_vc} and 
\ref{thm_uvc}, and illustrate their usefulness by disproving a conjecture of Xu.\footnote{Code 
that performs this calculation along with an implementation in SINGULAR 3.0
(http://www.singular.uni-kl.de) of the algorithms in this section can be found at
http://www.math.tamu.edu/$\sim$chillar/.} We also present some data to illustrate
their runtimes under different circumstances.

From Theorem \ref{thm_vc}, we have the following four methods for determining
$k$-colorability. They take as input a graph $G$ with vertices $V=\set{1,
\ldots, n}$ and edges $E$, and a positive integer $k$, and output \textsc{True}
if $G$ is $k$-colorable and otherwise \textsc{False}.

\vspace{0.5cm}
\begin{algorithmic}[1]
\Function{IsColorable}{$G$, $k$} \ [Theorem \ref{thm_vc} (2)] 
\State Compute a Gr\"obner basis $\mathcal{G}$ of $I_{G,k}$.
\State Compute the vector space dimension of $R/I_{G,k}$
over $\k$. \State \textbf{if} $\dim_\k R/I_{G,k} = 0$ \textbf{then return}
\textsc{False} \textbf{else return} \textsc{True}.
\EndFunction
\end{algorithmic}

\vspace{0.5cm}
\begin{algorithmic}[1]
\Function{IsColorable}{$G$, $k$} \ [Theorem \ref{thm_vc} (3)] 
\State Compute a Gr\"obner basis $\mathcal{G}$ of $I_{G,k}$.
\State Compute the normal form $\textup{nf}_\mathcal{G}(1)$ of
\Statex \hspace*{0.5cm}the constant polynomial 1 with respect to $\mathcal{G}$.
\State \textbf{if} $\textup{nf}_\mathcal{G}(1) = 0$ \textbf{then return}
\textsc{False} \textbf{else return} \textsc{True}.
\EndFunction
\end{algorithmic}

\vspace{0.5cm}
\begin{algorithmic}[1]
\Function{IsColorable}{$G$, $k$} \ [Theorem \ref{thm_vc} (4)] 
\State Set $\mathcal{G} := \set{x_i^k-1 : i \in V}$.
\State Compute the normal form $\textup{nf}_\mathcal{G}(f_G)$ of \label{nf2}
\Statex \hspace*{0.5cm}the graph polynomial $f_G$ with respect to $\mathcal{G}$.
\State \textbf{if} $\textup{nf}_\mathcal{G}(f_G) = 0$ \textbf{then return}
\textsc{False} \textbf{else return} \textsc{True}.
\EndFunction
\end{algorithmic}

\vspace{0.5cm}
\begin{algorithmic}[1]
\Function{IsColorable}{$G$, $k$} \ [Theorem \ref{thm_vc} (5)] 
\State Let $\mathcal{H}$ be the set of graphs with vertices $\set{1, \ldots, n}$
\Statex\hspace*{0.5cm}consisting of a clique of size $k+1$ and isolated vertices.
\State Set $\mathcal{G} := \set{f_H : H \in \mathcal{H}}$.
\State Compute the normal form $\textup{nf}_\mathcal{G}(f_G)$ of \label{nf3}
\Statex \hspace*{0.5cm}the graph polynomial $f_G$ with respect to $\mathcal{G}$.
\State \textbf{if} $\textup{nf}_\mathcal{G}(f_G) = 0$ \textbf{then return}
\textsc{False} \textbf{else return} \textsc{True}.
\EndFunction
\end{algorithmic}

\ \

From Theorem \ref{thm_uvc}, we have the following three methods for determining
unique $k$-colorability. They take as input a graph $G$ with vertices $V=\set{1, \ldots, n}$ and
edges $E$, and output \textsc{True}
if $G$ is uniquely $k$-colorable and otherwise \textsc{False}. 
Furthermore, the first two methods take as input a proper $k$-coloring $\nu$
of $G$ that uses all $k$ colors, while the last method requires a positive
integer $k$.

\vspace{0.5cm}
\begin{algorithmic}[1]
\Function{IsColorable}{$G$, $\nu$} \ [Theorem \ref{thm_uvc} (3)] 
\State Compute a Gr\"obner basis $\mathcal{G}$ of $I_{G,k}$.
\For{$i \in V$}
\State Compute the normal form $\textup{nf}_\mathcal{G}(g_i)$ of
\Statex \hspace*{1cm} the polynomial $g_i$ with respect to $\mathcal{G}$.
\State \textbf{if} $\textup{nf}_\mathcal{G}(g_i) \neq 0$ \textbf{then return}
\textsc{False}.
\EndFor
\State \Return{\textsc{True}}.
\EndFunction
\end{algorithmic}

\vspace{0.5cm}
\begin{algorithmic}[1]
\Function{IsColorable}{$G$, $\nu$} \ [Theorem \ref{thm_uvc} (4)] 
\State Compute a Gr\"obner basis $\mathcal{G}$ of $I_{n,k}:\<g_1, \ldots, g_n\>$.
\State Compute the normal form $\textup{nf}_\mathcal{G}(f_G)$ of \label{nf4}
\Statex \hspace*{0.5cm}the graph polynomial $f_G$ with respect to $\mathcal{G}$.
\State \textbf{if} $\textup{nf}_\mathcal{G}(f_G) = 0$ \textbf{then return}
\textsc{True} \textbf{else return} \textsc{False}.
\EndFunction
\end{algorithmic}

\vspace{0.5cm}
\begin{algorithmic}[1]
\Function{IsColorable}{$G$, $k$} \ [Theorem \ref{thm_uvc} (5)] 
\State Compute a Gr\"obner basis $\mathcal{G}$ of $I_{G,k}$.
\State Compute the vector space dimension of $R/I_{G,k}$
over $\k$.
\State \textbf{if} $\dim_\k R/I_{G,k} = k!$ \textbf{then return}
\textsc{True} \textbf{else return} \textsc{False}.
\EndFunction
\end{algorithmic}

\begin{rem}
It is possible to speed up the above algorithms dramatically
by doing some of the computations iteratively. First of all, step 2 of methods (2) and
(3) of Theorem \ref{thm_vc}, and methods (3) and (5) of Theorem \ref{thm_uvc}
should be replaced by the following code
\begin{algorithmic}[1]
\Statex
\State Set $I := I_{n,k}$.
\For{$\set{i,j} \in E$}
\State Compute a Gr\"obner basis $\mathcal{G}$ of
$I+ \< x_i^{k-1}+x_i^{k-2}x_j^{\phantom{1}} + \cdots + x_i^{\phantom{1}}x_j^{k-2}+x_j^{k-1} \>$.
\State Set $I := \<\mathcal{G}\>$.
\EndFor 
\end{algorithmic}

\ \

Secondly, the number of terms in the graph polynomial $f_G$ when fully expanded may be
very large. The computation of the normal form
$\textup{nf}_\mathcal{G}(f_G)$
of the graph polynomial $f_G$ in methods (4) and (5) of Theorem \ref{thm_vc}, and
method (4) of Theorem \ref{thm_uvc} should therefore be replaced by the
following code
\begin{algorithmic}[1]
\Statex
\State Set $f := 1$.
\For{$\set{i,j} \in E$ with $i<j$}
\State Compute the normal form $\textup{nf}_\mathcal{G}((x_i-x_j)f)$ of
\Statex \hspace*{0.5cm}$(x_i-x_j)f$ with respect to $\mathcal{G}$, and set $f :=
\textup{nf}_\mathcal{G}((x_i-x_j)f)$.
\EndFor
\end{algorithmic}
\end{rem}

\begin{figure}[!htbp]
\begin{center}
\includegraphics[scale=1.0]{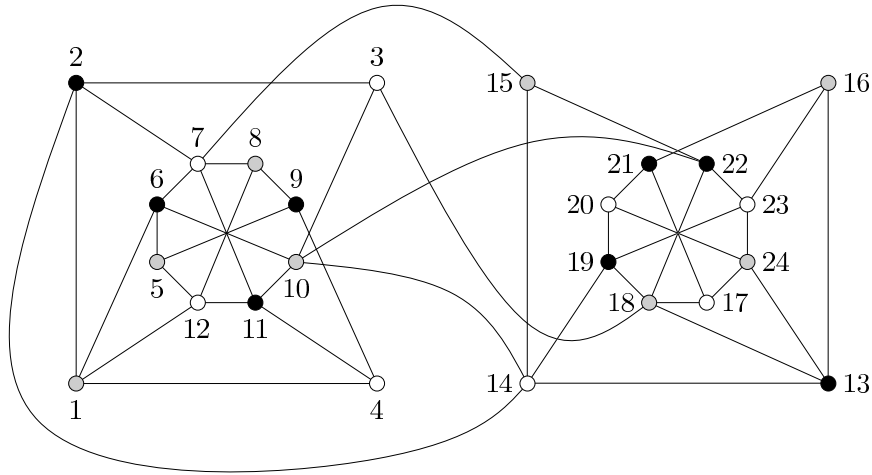}
\end{center}
\caption{A counterexample to Xu's conjecture \cite{akbari1}.}
\label{fig.xucounterex}
\end{figure}

In \cite{xu}, Xu showed that if $G$ is a uniquely $k$-colorable graph with
$|V|= n$ and $|E| = m$, then $m \geq (k-1)n-{k \choose 2}$, and
this bound is best possible.  He went on to conjecture that if $G$ is
uniquely $k$-colorable with $|V|= n$ and $|E| = (k-1)n-{k \choose 2}$, then $G$
contains a $k$-clique.  In \cite{akbari1}, this conjecture was
shown to be false for $k=3$ and $|V| = 24$ using the graph in Figure \ref{fig.xucounterex}; 
however, the proof is somewhat complicated.  We verified that this graph is
indeed a counterexample to Xu's conjecture using several of the above mentioned
methods. The fastest verification requires less than two seconds of processor time on a laptop PC with a $1.5$ GHz Intel Pentium M processor 
and $1.5$ GB of memory.  The code can be downloaded from the link
at the beginning of this section.  The speed of these calculations should
make the testing of conjectures for uniquely colorable graphs 
a more tractable enterprise.

Below are the runtimes for the graphs in Figures 1 and 2. The term orders used are
given in the notation of the computational algebra program Singular: 
lp is the lexicographical ordering, Dp
is the degree lexicographical ordering, and dp is the degree reverse
lexicographical ordering. That the computation did not finish within 10 minutes
is denoted by ``$> 600$'', while ``--'' means that the computation ran out of
memory.
 
 \ \

\begin{center}
\begin{tabular}{|c|ccc|ccc|}
\hline
Characteristic &  & 0 &  &  & 2 &  \\
\hline
Term order & lp & Dp & dp & lp & Dp & dp \\
\hline
Theorem \ref{thm_vc} (2) & 3.28 & 2.29 & 1.24 & 2.02 & 1.56 & 0.81 \\
Theorem \ref{thm_vc} (3) & 3.30 & 2.42 & 1.25 & 2.15 & 1.60 & 0.94 \\
Theorem \ref{thm_vc} (4) & 1.86 & $> 600$ & $> 600$ & 1.08 & 448.28 & 324.89 \\
Theorem \ref{thm_vc} (5) & $> 600$ & $> 600$ & $> 600$ & $> 600$ & $> 600$ & $> 600$ \\
\hline
Theorem \ref{thm_uvc} (3) & 3.53 & 2.54 & 1.43 & 2.23 & 1.72 & 1.03 \\
Theorem \ref{thm_uvc} (4) & $> 600$ & $> 600$ & $> 600$ & $> 600$ & $> 600$ & $> 600$ \\
Theorem \ref{thm_uvc} (5) & 3.30 & 2.28 & 1.24 & 2.03 & 1.54 & 0.82 \\
\hline
\end{tabular}

Runtimes in seconds for the graph in Figure 1.
\end{center}

\ \

\begin{center}
\begin{tabular}{|c|ccc|ccc|}
\hline
Characteristic &  & 0 &  &  & 2 &  \\
\hline
Term order & lp & Dp & dp & lp & Dp & dp \\
\hline
Theorem \ref{thm_vc} (2) & 596.89 & 33.32 & 2.91 & 144.05 & 12.45 & 1.64 \\
Theorem \ref{thm_vc} (3) & 598.25 & 33.47 & 2.87 & 144.60 & 12.44 & 1.81 \\
Theorem \ref{thm_vc} (4) & -- & $> 600$ & $> 600$ & -- & $> 600$ & $> 600$ \\
Theorem \ref{thm_vc} (5) & $> 600$ & $> 600$ & $> 600$ & $> 600$ & $> 600$ & $> 600$ \\
\hline
Theorem \ref{thm_uvc} (3) & 597.44 & 34.89 & 4.29 & 145.81 & 13.55 & 3.02 \\
Theorem \ref{thm_uvc} (4) & -- & -- & -- & -- & -- & -- \\
Theorem \ref{thm_uvc} (5) & 595.97 & 33.46 & 2.94 & 145.02 & 12.34 & 1.64 \\
\hline
\end{tabular}

Runtimes in seconds for the graph in Figure 2.
\end{center}

\ \

Another way one can prove that a graph is uniquely $k$-colorable
is by computing the chromatic polynomial and testing if it equals $k!$ when
evaluated at $k$. This is possible for the graph in Figure 1.  Maple
reports that it has chromatic polynomial
\begin{multline*}
x(x-2)(x-1)(x^9-20x^8+191x^7-1145x^6+4742x^5\\
-14028x^4+29523x^3-42427x^2+37591x-15563).
\end{multline*}
When evaluated at $x=3$ we get the expected
result $6 = 3!$. Computing the above chromatic polynomial took 94.83 seconds.
Maple, on the other hand, was not able to
compute the chromatic polynomial of the graph in Figure 2 within 10 hours. 


\section{Acknowledgments}

We would like to thank the anonymous referees for valuable comments that
greatly improved this work.  We also thank Keith Briggs for pointing out an error
in our drawing of the graph appearing in Figure 2.

\end{document}